\documentstyle[amssymb,12pt]{article}

\newtheorem{theorem}{\sc Theorem}

\newtheorem{definition}[theorem]{\sc Definition}

\begin{document}

\author{Paul C. Eklof\thanks{%
Thanks to Rutgers University for its support of this research through its
funding of the authors' visits to Rutgers.} \\
University of California, Irvine \and Saharon Shelah\thanks{%
Partially supported by Basic Research Fund, Israeli Academy of Sciences.
Pub. No. 559a} \\
Hebrew University and Rutgers University}
\title{New non-free Whitehead groups \\
(corrected version)}
\date{}
\maketitle

\begin{abstract}
We answer an open problem in \cite{EM} by showing that it is consistent that
there is a strongly $\aleph _{1}$-free $\aleph _{1}$-coseparable group of
cardinality $\aleph _{1}$ which is not $\aleph _{1}$-separable.
\end{abstract}

\section{INTRODUCTION}

An abelian group $A$ is called a Whitehead group, or $W$-group for short, if 
$\hbox{Ext}(A,{\Bbb {Z}})=0$. For historical reasons, $A$ is called an $%
\aleph _{1}$-coseparable group if $\hbox{Ext}(A,{\Bbb {Z}}^{(\omega )})=0$%
, but for convenience we shall use non-standard terminology and say $A$ is a 
$W_{\omega }$-group when $A$ is $\aleph _{1}$-coseparable. Obviously a $%
W_{\omega }$-group is a $W$-group. In 1973--75, the second author proved
that it is consistent with ZFC + GCH that every $W$-group is free and
consistent with ZFC that there are non-free $W_{\omega }$-groups of
cardinality $\aleph _{1}$ (\cite{Sh44}, \cite{Sh75}); he later showed that
it is consistent with ZFC + GCH that there are non-free $W_{\omega }$-groups
of cardinality $\aleph _{1}$ (\cite{Sh77}, \cite{Sh80}). Before 1973 it was
known (in ZFC) that every $W$-group is $\aleph _{1}$-free, separable, and
slender, and assuming CH, every $W$-group is strongly $\aleph _{1}$-free.
(See, for example, \cite[pp. 178--180]{F}.) These turned out, by the results
of the second author, to be essentially all that could be proved without
additional set-theoretic hypotheses.

However, new questions of what could be proved in ZFC arose, inspired by the
consistency results and their proofs. One of the most intriguing was:

\begin{quote}
(0.1) Does every strongly $\aleph _1$-free $W_\omega $-group of cardinality $%
\aleph _1$ satisfy the stronger property that it is $\aleph _1$-separable?
\end{quote}

\noindent (See \cite[p. 454, Problem 5]{EM}). As we shall explain below, not
only was the answer to this question affirmative in every known model of
ZFC, but the nature of the known constructions of non-free Whitehead groups
was such as to lead to the suspicion that the answer might be affirmative
(provably in ZFC). However, in this paper we show that it is consistent that
the answer is negative.

First we recall the key definitions. An abelian group $A$ is $\aleph _1${\it %
-free} if every countable subgroup of $A$ is free;{\it \ } $A$ is {\it %
strongly }$\aleph _1${\it -free} if every countable subset is contained in a
countable free subgroup $B$ such that $A/B$ is $\aleph _1$-free. $A$ is $%
\aleph _1${\it -separable} if every countable subset is contained in a
countable free subgroup $B$ which is a direct summand of $A$; so an $\aleph
_1$-separable group is strongly $\aleph _1$-free. It is a consequence of CH
(or even of $2^{\aleph _0}<2^{\aleph _1}$) that there are strongly $\aleph
_1 $-free groups of cardinality $\aleph _1$ which are not $\aleph _1$%
-separable (see \cite{Sh81}). However, the existence of such groups is not
settled by the hypothesis $2^{\aleph _0}=2^{\aleph _1}$; specifically, in a
model of MA + $\neg $CH every strongly $\aleph _1$-free group of cardinality 
$\aleph _1$ {\it is} $\aleph _1$-separable; but the methods of \cite{MS}
show that it is consistent with $2^{\aleph _0}=2^{\aleph _1}$ that there are
strongly $\aleph _1$-free groups of cardinality $\aleph _1$ which are {\it %
not} $\aleph _1$-separable.

Now suppose $A$ is strongly $\aleph _{1}$-free and is a $W_{\omega }$-group.
Consider a countable subgroup $B$ of $A$ such that $A/B$ is $\aleph _{1}$%
-free. We have a short exact sequence 
\[
0\rightarrow B\rightarrow A\rightarrow A/B\rightarrow 0 
\]
where the map of $B$ into $A$ is inclusion. Since $B$ is a free group of
countable rank, if we knew that $A/B$ were a $W_{\omega }$-group, then we
would have $\hbox{Ext}(A/B,{\Bbb {Z}}^{(\omega )})=\hbox{Ext}(A/B,B)=0$
and we could conclude that this sequence splits and hence $B$ is a direct
summand of $A$. In every previously known model where there are non-free $%
W_{\omega }$-groups, the construction of a $W_{\omega }$-group $A$ is such
that $A/B$ shares the properties of $A$ closely enough that $A/B$ is also a $%
W_{\omega }$-group --- when $B$ is a countable subgroup such that $A/B$ is $%
\aleph _{1}$-free. (For example, if $A$ is constructed as in 
\cite[Prop. XII.3.6(iii), p. 371]{EM}, using a ladder system with a
uniformization property, then $A/B$ shares the same properties, because it
is constructed using essentially the same ladder system.) Thus in these
models the answer to (0.1) is affirmative. This motivates question (0.1) as
well as the related question

\begin{quote}
(0.2) If a group $A$ of cardinality $\aleph _1$ is strongly $\aleph _1$-free
and a $W_\omega $-group, and $B$ is a countable subgroup of $A$ such that $%
A/B$ is $\aleph _1$-free, is $A/B$ a $W_\omega $-group?
\end{quote}

\noindent By what we have just remarked, a positive answer to (0.2) implies
a positive answer to (0.1). The converse holds as well: if $A$ and $B$ are
as in the hypotheses of (0.2) and $A$ is $\aleph _{1}$-separable, $A=F\oplus
A^{\prime }$ where $F$ is countable and contains $B$; then $A/B$ is a $%
W_{\omega }$-group because $A/B=F/B\oplus A^{\prime }$ and $F/B$ is free by
hypothesis on $B$.

We shall give a model of ZFC +$\lnot $CH where the answer to (0.1) and (0.2)
is negative.\footnote{%
We do not know if Theorem 8 of the original paper is correct, or if the
answer to question (0.3) is ``no''.}

\section{THE PROOF}

Our main theorem is:

\begin{theorem}
There is a strongly $\aleph _{1}$-free $W_{\omega }$-group $A$ of
cardinality $\aleph _{1}$ with a countable subgroup $B$ of $A$ such that $%
A/B $ is $\aleph _{1}$-free but $B$ is not a direct summand of $A$.
\end{theorem}

Throughout, $E$ will be a stationary subset of $\omega _{1}$ consisting of
limit ordinals, with (for technical reasons) $\omega \notin E$. We begin
with a general construction of a group. Let $\pi _{n}$ be the $n$th prime.

\begin{definition}
\label{A} For each $\delta \in E$ let $\eta _{\delta }$ be a ladder on $%
\delta $, that is, a strictly increasing function $\eta _{\delta }:\omega
\rightarrow \delta $ whose range approaches $\delta $. Let $\varphi $ be a
function from $E\times \omega $ to $\omega $. Let $F$ be the free abelian
group with basis $\{x_{\nu }\colon \nu \in \omega _{1}\}\cup \{z_{\delta ,n}%
\colon \delta \in E,n\in \omega \}$ and let $K$ be the subgroup of $F$
generated by $\{w_{\delta ,n}:\delta \in E,n\in \omega \}$ where 
\begin{equation}
\hbox{ }w_{\delta ,n}=\pi _{n}z_{\delta ,n+1}-z_{\delta ,0}-x_{\eta _{\delta
}(n)}-x_{\varphi (\delta ,n)}.  \label{1.1}
\end{equation}
Let $A=F/K$.
\end{definition}

Clearly $A$ is an abelian group of cardinality $\aleph _{1}$. Notice that
because the right-hand side of (\ref{1.1}) is 0 in $A$, we have for each $%
\delta \in E$ and $n\in \omega $ the following relation in $A$: 
\begin{equation}
\pi _{n}z_{_{\delta ,n+1}}=z_{\delta ,0}+x_{\eta _{\delta }(n)}+x_{\varphi
(\delta ,n)}  \label{1.2}
\end{equation}
where, in an abuse of notation, we write, for example, $z_{_{\delta ,n+1}}$
instead of $z_{_{\delta ,n+1}}+K$. If we let 
\begin{equation}
\hbox{ }A_{\alpha }=\left\langle \{x_{\nu }:\nu <\alpha \}\cup \{z_{\delta
,n}:\delta \in E\cap \alpha \hbox{, }n\in \omega \}\right\rangle \hbox{.}
\label{1.4}
\end{equation}
for each $\alpha <\omega _{1}$, then for each $\delta \in E$, $z_{\delta
,0}+A_{\delta }$ is non-zero and divisible in $A_{\delta +1}/A_{\delta }$ by 
$\pi _{n}$ for all $n\in \omega $. Thus $A_{\delta +1}/A_{\delta }$ is not
free and hence $A$ is not free. (In fact $\Gamma (A)\supseteq \tilde{E}$;
see \cite[pp. 85f]{EM}.) Moreover, $A$ is strongly $\aleph _{1}$-free; in
fact, for every $\alpha <\omega _{1}$, using Pontryagin's Criterion we can
show that $A/A_{\alpha }$ is $\aleph _{1}$-free whenever $\alpha \notin E$.

We now define the model of ZFC where $A$ is defined and has the desired
properties. We begin with a model $V$ of ZFC where GCH holds, choose $E\in V$%
, and define the group $A$ in a generic extension $V^{Q_{0}}$ using generic
ladders $\eta _{\delta }$, and generic $\varphi $. Specifically:

\begin{definition}
\label{Q0}Let $Q_{0}$ be the set of all finite functions $q$ such that $%
\hbox{dom}(q)$ is a finite subset of $E$ and for all $\gamma \in \hbox{%
dom}(q)$, $q(\gamma )$ is a pair $(\eta _{\gamma }^{q},\varphi _{\gamma
}^{q})$ where for some $r_{\gamma }^{q}\in \omega $:

\begin{itemize}
\item  $\eta _{\gamma }^{q}$ is a strictly increasing function$:r_{\gamma
}^{q}\rightarrow \gamma $;

\item  $\varphi _{\gamma }^{q}:\{\gamma \}\times r_{\gamma }^{q}\rightarrow
\omega $.
\end{itemize}
\end{definition}

Clearly $Q_{0}$ is c.c.c. We now do an iterated forcing to make $A$ a $%
W_{\omega }$-group. We begin by defining the basic forcing that we will
iterate.

\begin{definition}
Given a homomorphism $\psi :K\rightarrow {\Bbb Z}^{(\omega )}$, let $Q_{\psi
}$ be the poset of all finite functions $q$ into ${\Bbb Z}^{(\omega )}$
satisfying:

There are $\delta _{0}<\delta _{1}<...<\delta _{m}$ in $E$ and $\{r_{\ell
}:\ell \leq m\}\subseteq \omega $ such that $\hbox{dom}(q)=$%
\[
\{z_{\delta _{\ell },n}:\ell \leq m,n\leq r_{\ell }\}\cup \{x_{\nu }:\nu \in
I_{q}\} 
\]
where $I_{q}\subset \omega _{1}$ is finite and is such that for all $\ell
\leq m$

\begin{equation}
\eta _{\delta _{\ell }}(n)\in I_{q}\Leftrightarrow n<r_{\ell }  \label{5}
\end{equation}

\noindent and for all $\ell \leq m$ and $n<r_{\ell }$,

\begin{equation}
\psi (w_{\delta _{\ell },n})=\hbox{ }\pi _{n}q(z_{\delta _{\ell
},n+1})-q(z_{\delta _{\ell },0})-q(x_{\eta _{\delta _{\ell
}}(n)})-q(x_{\varphi (\delta ,n)}).  \label{6}
\end{equation}

Moreover, we require of $q$ that for all $\ell \neq j$ in $\{0,...,m\}$, 
\begin{equation}
\hbox{ }\eta _{\delta _{j}}(k)\neq \eta _{\delta _{\ell }}(i)\hbox{ for all }%
k\geq r_{j}\hbox{ and }i\in \omega .  \label{7}
\end{equation}
\end{definition}

We will denote $\{\delta _{0},...,\delta _{m}\}$ by $\hbox{cont}(q)$ and $%
r_{\ell }$ by $\hbox{num}(q,\delta _{\ell })$. The partial ordering on $%
Q_{\psi }$ is inclusion. Standard methods prove that $Q_{\psi }$ is c.c.c.

\smallskip\ 

Let $P=\left\langle P_{i},\dot{Q}_{i}:0\leq i<\omega _{2}\right\rangle $ be
a finite support iteration of length $\omega _{2}$ so that for every $i\geq
1 $ $\Vdash _{P_{i}}\dot{Q}_{i}=Q_{\dot{\psi}_{i}}$ where $\Vdash _{P_{i}}%
\dot{\psi}_{i}$ is a homomorphism$:K\rightarrow {\Bbb Z}^{(\omega )}$, where
the enumeration of names $\{\dot{\psi}_{i}:1\leq i<\omega _{2}\}$ is chosen
so that if $G$ is $P$-generic and $\psi \in V[G]$ is a homomorphism$%
:K\rightarrow {\Bbb Z}^{(\omega )}$, then for some $i\geq 1$, $\dot{\psi}%
_{i} $ is a name for $\psi $ in $V^{P_{i}}$. Then $P$ is c.c.c. and in $V[G]$
every homomorphism from $K$ to ${\Bbb Z}^{(\omega )}$ extends to one from $F$
to ${\Bbb Z}^{(\omega )}$. This means that $\hbox{Ext}(A,{\Bbb Z}%
^{(\omega )})=0$, that is, $A$ is a $W_{\omega }$-group (see, for example, 
\cite[p.8]{EM}).

Let $B=A_{\omega }$, i.e., the subgroup of $A$ generated by $\{x_{\ell
}:\ell \in \omega \}$. It is easy to check that $A/B$ is $\aleph _{1}$-free.
\smallskip Now, aiming for a contradiction, suppose that in $V[G]$ there\ is
a projection $h:A\rightarrow B$ (i.e., $h\upharpoonright B$ is the
identity). Then there is a condition $p_{o}\in G$ such that 
\[
p_{o}\Vdash \hbox{``}\dot{h}:A\rightarrow B\hbox{ is a projection''} 
\]
where $\dot{h}$ is a name for $h$.

For each ordinal $\xi \in \omega _{1}-E$, choose a condition $p_{\xi }\geq
p_{o}$ such that there is a $y_{\xi }\in V$ such that 
\[
p_{\xi }\Vdash \dot{h}(x_{\xi +1})=y_{\xi }. 
\]
(That is, $y_{\xi }$ is an element of $B\cap V$, and not just a name.)

We can assume that

\begin{quote}
\bigskip \ \ ($\dag $) $0\in \hbox{dom}(p_{\xi })$; for each $j\in 
\hbox{dom}(p_{\xi })$, $p_{\xi }(j)$ is a function in $V$ and not just a
name; $r_{\gamma }^{p_{\xi }(0)}$($=r_{\xi }$) is independent of $\gamma \in 
\hbox{dom}(p_{\xi }(0))$; if $j\in \hbox{dom}(p_{\xi })\setminus \{0\} 
$, $\gamma \in \hbox{cont}(p_{\xi }(j))$ implies $\gamma \in \hbox{dom}%
(p_{\xi }(0))$ and $\hbox{num}(p_{\xi }(j),\gamma )$ ($=r_{\xi
,j}^{\prime }$) is $\leq r_{\xi }$ and independent of $\gamma $. Moreover,
if $\gamma >\xi $, then $\eta _{\gamma }(r_{\xi ,j}^{\prime }-1)>\xi $.
\end{quote}

When we say that ``$\nu $ occurs in $p$'' we mean that $\nu \in \hbox{dom}%
(p)\cup \hbox{dom}(p(0))\cup \bigcup \{\{\eta _{\gamma }(n),\varphi
(\delta ,n)\}:\gamma \in \hbox{dom}(p(0)$, $n<r_{\gamma }^{p(0)}\}$ or $%
x_{\nu }$ belongs to the domain of some $p(j)$.

Without loss of generality we can assume (passing to a subset $S\subseteq
\omega _{1}-E$) by Fodor's Lemma and the $\Delta $-system lemma that

\begin{quote}
($\dag \dag $) $\{\hbox{dom}(p_{\xi }):\xi \in S\}$ forms a $\Delta $%
-system, whose root we denote $C$ (i.e., $\hbox{dom}(p_{\xi _{1}})\cap 
\hbox{dom}(p_{\xi _{2}})=C=\{0,\mu _{1},...,\mu _{d}\}$ for all $\xi
_{1}\neq \xi _{2}$ in $S$); $r_{\xi }$ ($=r$) and $r_{\xi ,j}^{\prime }$ ($%
=r_{j}^{\prime }$) are independent of $\xi $; $\hbox{dom}(p_{\xi
}(0))\cap \xi $ is independent of $\xi $; there is $m$ such that for all $%
\xi \in S$, $\hbox{dom}(p_{\xi }(0))-\xi =\{\gamma _{\xi ,0}<...<\gamma
_{\xi ,m}\}$and for each $\ell \leq m$ there is $t_{\ell }\leq r-1$ such
that $\eta _{\gamma _{\xi ,\ell }}(t_{\ell }-1)<\xi \leq \eta _{\gamma _{\xi
,\ell }}(t_{\ell })$ and for $i\leq t_{\ell }-1$, $\eta _{\gamma _{\xi ,\ell
}}(i)$ is independent of $\xi $. Moreover, for every $j\in C$, $\{\hbox{%
dom}(p_{\xi }(j)):\xi \in S\}$ forms a $\Delta $-system and for all $\xi
_{1}\neq \xi _{2}$ in $S$, $p_{\xi _{1}}(j)$ and $p_{\xi _{2}}(j)$ agree on $%
\hbox{dom}(p_{\xi _{1}}(j))\cap \hbox{dom}(p_{\xi _{2}}(j))$.
\end{quote}

Let $p^{*}$ denote the ``heart'' of the $\Delta $-system; that is, $\hbox{%
dom}(p^{*})=C$ and for all $\mu \in C$, $\hbox{dom}(p^{*}(\mu ))=\hbox{%
dom}(p_{\xi _{1}}(\mu ))\cap \hbox{dom}(p_{\xi _{2}}(\mu ))$ (= $C_{\mu }$%
, say) for $\xi _{1}\neq \xi _{2}\in S$; and $p^{*}(\mu )\upharpoonright
C_{\mu }=p_{\xi _{1}}(\mu )\upharpoonright C_{\mu }$.

We can assume that every ordinal which occurs in $p^{*}$ is $<\xi $ for
every $\xi \in S$, and that for every $\xi _{1}<\xi _{2}$ in $S$, every $\nu 
$ which occurs in $\xi _{1}$ is $<\xi _{2}$. We can find $\delta \in E$
which is the limit of a strictly increasing sequence $\{\xi _{\ell }:\ell
\in \omega \}\subseteq S$. Notice that no ordinal $\geq $ $\delta $ occurs
in any $p_{\xi _{\ell }}$. There is a condition $\tilde{p}\geq p^{*}$ which
forces a value to $\dot{h}(z_{\delta ,0})$, i.e., there is $y\in V$ such
that $\tilde{p}\Vdash \dot{h}(z_{\delta ,0})=y$. We can assume that $\tilde{p%
}$ is as in ($\dag $); let $r^{*}=r_{\delta }^{\tilde{p}(0)}$. Fix $n\in
\omega $ sufficiently large so that every ordinal which occurs in $\tilde{p}$
but not in $p^{*}$ does not occur in $p_{\xi _{n}}$. Fix $k\in \omega $ such
that $k>\sup \{j\in \omega :j$ occurs in $\tilde{p}$ or $p_{\xi _{n}}\}$ and
such that $y,y_{\xi _{n}}\in \{x_{j}:j<k\}$. There is a condition $q_{0}$
extending $\tilde{p}(0)$ and $p_{\xi _{n}}(0)$ such that $q_{0}$ forces 
\[
\eta _{\delta }(r^{*})=\xi _{n}+1\hbox{; }\eta _{\delta }(r^{*}+1)=\xi _{n+1}%
\hbox{; and }\varphi (\delta ,r^{*})=k\hbox{.} 
\]
Then there is a condition $q\in P$ extending $\tilde{p}$ and $p_{\xi _{n}}$
such that $q(0)\geq q_{0}$. (The only possible difficulty in defining $q(j)$
for $j>0$ is in defining $q(j)(z_{\delta ,r^{*}+1})$, $q(j)(\eta _{\delta
}(r^{*}))$ and $q(j)(x_{\varphi (\delta ,r^{*})})$ to satisfy (\ref{6}), but
this can be done even though $q(j)(\eta _{\delta }(r^{*}))=q(j)(x_{\xi
_{n}+1})$ may be determined by $p_{\xi _{n}}(j)$, because $x_{k}$ is new.)

Now consider a generic extension $V[G^{\prime }]$ such that $q\in G^{\prime
} $. In this model, by (\ref{1.2}), we have that in $A$, $\pi _{n}$ divides $%
z_{\delta ,0}-x_{\eta _{\delta }(r^{*})}-x_{\varphi (\delta
,r^{*})}=z_{\delta ,0}-x_{\xi _{n}+1}-x_{k}$. Hence in $B$, $\pi _{n}$
divides $h(z_{\delta ,0}-x_{\xi _{n}+1}-x_{k})=y-y_{\xi _{n}}-x_{k}$. But
this is impossible by choice of $k$.

\end{document}